\documentclass[sn-mathphys,Numbered]{sn-jnl}

\usepackage{graphicx}%
\usepackage{multirow}%
\usepackage{amsmath,amssymb,amsfonts}%
\usepackage{amsthm}%
\usepackage{mathrsfs}%
\usepackage[title]{appendix}%
\usepackage{xcolor}%
\usepackage{textcomp}%
\usepackage{manyfoot}%
\usepackage{booktabs}%
\usepackage{indentfirst}
\usepackage{algpseudocode}%
\usepackage{listings}%
\usepackage{hyperref}
\usepackage[ruled,vlined]{algorithm2e}
\usepackage{setspace}

\begin{document}

\title[Article Title]{An adaptive linearized alternating direction multiplier method with a relaxation step  for  convex programming}


\author{\fnm{Boran Wang}}\email{1932210936@qq.com}

\affil{\orgdiv{College of Science}, \orgname{Minzu University of China},  \city{Beijing}, \postcode{100081}, \state{China}}


\abstract{Alternating direction multiplication is a powerful technique for solving convex optimisation problems. When challenging subproblems are encountered in the real world, it is useful to solve them by introducing neighbourhood terms. When the neighbourhood matrix is positive definite, the algorithm converges but at the same time makes the iteration step small. Recent studies have revealed the potential non-positive definiteness of the neighbourhood matrix. In this paper, we present an adaptive linearized alternating direction multiplier method with a relaxation step, combining the relaxation step with an adaptive technique. The novelty of the method is to use the information of the current iteration point to dynamically select the neighbourhood matrix, increase the iteration step size, and speed up the convergence of the algorithm.We prove the global convergence of the algorithm theoretically and illustrate the effectiveness of the algorithm using numerical experiments.}

\keywords{variational inequality, an adaptive linearized ADMM , a relaxation step , global convergence.}



\maketitle

\section{Introduction}\label{sec1}
 
This paper considers the following separable convex optimization problem with linear constraints:
\begin{equation}\label{eqn:1} 
\min\left \{f\left ( x \right ) +  g\left ( y \right ) | Ax+By=b,x\in X,y\in Y  \right \} 
\end{equation}
where $f:R^{n_{1} } \to R$ and $g:R^{n_{2} } \to R$ are convex functions (but not necessarily smooth), $A\in R^{m\times n_{1} } ,B\in R^{m\times n_{2} } ,b\in R^{m},X\subseteq R^{n_{1} }$  and $Y\subseteq R^{n_{2} }$ are closed convex sets. Problem \eqref{eqn:1} has recently been widely used in several areas such as image processing\cite{1,2}, statistical learning\cite{3} and communication networks\cite{4}.There exist numerous effective methods for addressing problem \eqref{eqn:1}, including the primal-dual method,  the augmented Lagrange method, and the alternating direction multiplier method. Notably, the alternating direction multiplier method has garnered significant attention in recent years due to its simplicity and efficiency. Many scholars have introduced several effective variations of this algorithm tailored to practical challenges.

The alternate direction multiplier approach, which was first put forth by Gabay and Mercier\cite{5} and Glowinski and Marrocco\cite{6}, is widely recognized as an effective way to solve \eqref{eqn:1}. For the progress of research on the alternating direction multiplier method, we refer to\cite{7}.The method possesses fast convergence speed and strong convergence performance, making it a crucial tool for addressing convex optimization problems with divisibility. It can decompose a large problem into two or more small-scale subproblems and then iteratively solve each subproblem to enhance solution speed. Its direct, adaptable, and practical nature is particularly effective in tackling convex optimization problems with separability. The alternating direction multiplier method can be described as follows:

\begin{equation}\label{eqn:2} 
\left\{
\begin{aligned}
x^{k+1} &=arg\min_{x}\left \{ L_{\beta } \left ( x,y^{k},\lambda ^{k}   \right )|x\in X  \right \}
 \\y^{k+1} &=arg\min_{y}\left \{L_{\beta } \left ( x^{k+1} ,y,\lambda ^{k}   \right )| y\in Y \right \}
 \\\lambda ^{k+1}& =\lambda ^{k} -\beta \left ( Ax^{k+1} +By^{k+1} -b \right ) 
\end{aligned}
\right.
\end{equation}
where $L_{\beta }\left ( x,y,\lambda  \right )$ is the augmented Lagrangian  function of \eqref{eqn:1}
\\$$L_{\beta }\left ( x,y,\lambda  \right )=f\left ( x \right) +g\left ( y \right )-\lambda ^{T}\left ( Ax+By-b \right )+\frac{\beta }{2}\left \| Ax+By-b \right \|_{2}^{2}.$$\\
where $\lambda \in R^{m}$ represents the Lagrange multiplier, and $ \beta > 0$ represents the penalty parameter. For simplicity, the penalty parameter $ \beta $ is fixed in our discussion.

A significant focus in the diverse research literature concerning the alternating direction multiplier method is the exploration of efficient strategies for solving its subproblems. In certain scenarios, the functions $f\left ( x \right )$ and $g\left ( y \right )$, or the coefficient matrices $A$ and $B$ may exhibit unique properties and structures. Leveraging these characteristics, researchers can extend the framework presented in \eqref{eqn:2} to devise algorithms tailored to specific applications while upholding theoretical convergence guarantees. This principle underscores the critical role of effectively applying the alternating direction multiplier method across various fields and applications. To delve deeper into this concept, let's examine the $y$-subproblem outlined in \eqref{eqn:2}. As previously noted, the  function $g(y)$, the  matrice $B$, and the set $Y$ are pivotal in addressing the $y$-subproblems. Typically, when these the function, matrice, and set are in a generic form, the solution to \eqref{eqn:2} can be obtained using straightforward iterative techniques,we refer to \cite{8,9,10}. However, in practical scenarios, the functions, matrices, and sets $Y$ associated with the $y$-subproblem may possess distinct characteristics, necessitating more efficient solution methods. When the matrix B is not an identity matrix, the $y$-subproblem \eqref{eqn:2} can be reformulated as follows: 
\begin{equation}\label{eqn:3} 
 y^{k+1}=arg\min_{y}\left \{ g\left ( y \right ) +  \frac{\beta }{2} \left \| By+\left ( Ax^{k+1}-b-\frac{1}{\beta } \lambda ^{k}   \right )  \right \|^{2}  | y\in Y \right \}   
\end{equation}

After linearizing the quadratic term $\left \| By+ \left ( Ax^{k+1} -b-\frac{1}{\beta }\lambda ^{k} \right ) \right \|^{2}$ in \eqref{eqn:3} , it is obtained:
\begin{equation}\label{eqn:4}
y^{k+1}=arg\min_{y}\left \{ g\left ( y \right )+  \frac{r}{2} \left \| y-\left ( y^{k}+  \frac{1}{r}q^{k}    \right )  \right \| ^{2} |y\in  Y    \right \}.    
\end{equation}
where
\begin{equation}\label{eqn:5}
q^{k}=B^{T}\left ( \lambda ^{k}-\beta \left ( Ax^{k+1}+By^{k}-b   \right )   \right ).    
\end{equation}
where $ r> 0 $ is a constant, and thus we obtain the linearized alternating direction multiplier method\cite{11}:
\begin{equation}\label{eqn:6}
\left\{
\begin{aligned}
x^{k+1} &=arg\min_{x}\left \{ L_{\beta } \left ( x,y^{k},\lambda ^{k}   \right )|x\in X  \right \}
 \\y^{k+1} &=arg\min_{y}\left \{L_{\beta } \left ( x^{k+1} ,y,\lambda ^{k}   \right )+\frac{r}{2}\left \| y-\left (  y^{k}+\frac{1}{r}q^{k}   \right )  \right \|^{2}   | y\in Y \right \}
 \\\lambda ^{k+1}& =\lambda ^{k} -\beta \left ( Ax^{k+1} +By^{k+1} -b \right ) 
\end{aligned}
\right.
\end{equation}
Linearized alternating direction multiplier method is widely used with compressed perception\cite{12}, and image processing\cite{13,14}, etc. . 

A more general  alternating direction multiplier method\cite{15} can be written as follows:
\begin{equation}\label{eqn:7}
\left\{
\begin{aligned}
x^{k+1} &=arg\min_{x}\left \{ L_{\beta } \left ( x,y^{k},\lambda ^{k}   \right )|x\in X  \right \}
 \\y^{k+1} &=arg\min_{y}\left \{L_{\beta } \left ( x^{k+1} ,y,\lambda ^{k}   \right )+\frac{1}{2}\left \| y-y^{k}    \right \|_{D} ^{2}   | y\in Y \right \}
 \\\lambda ^{k+1}& =\lambda ^{k} -\beta \left ( Ax^{k+1} +By^{k+1} -b \right ) 
\end{aligned}
\right.
\end{equation}
where $D\in R^{n_{2} \times n_{2} }$ is positive definite, and thus, the linearized alternating direction multiplier method can be viewed as a special case of the general  alternating direction multiplier method, where the positive definite proximal term.
\begin{equation}
D=rI_{n_{2} }-\beta B^{T}B\quad and \quad r> \beta \left \| B^{T}B  \right \|    
\end{equation}
Since in many practical applications, only one subproblem in \eqref{eqn:6} needs to be linearized, in this paper, we only need to consider the case of linearizing the $y$-subproblem in \eqref{eqn:7}.

In order to improve the applicability of the linearized alternating direction multiplier method, a number of improved linearized alternating direction multiplier methods\cite{16,17,18,19,20}have been proposed by some scholars.

Literature\cite{21} proposes an indefinite proximity alternating direction multiplier method with the following iteration format:
\begin{equation}
\left\{
\begin{aligned}
x^{k+1} &=arg\min_{x}\left \{ L_{\beta } \left ( x,y^{k},\lambda ^{k}   \right )|x\in X  \right \}
 \\y^{k+1} &=arg\min_{y}\left \{L_{\beta } \left ( x^{k+1} ,y,\lambda ^{k}   \right )+\frac{1}{2}\left \| y-y^{k}    \right \|_{D_{0} } ^{2}   | y\in Y \right \}
 \\\lambda ^{k+1}& =\lambda ^{k} -\beta \left ( Ax^{k+1} +By^{k+1} -b \right ) 
\end{aligned}
\right.
\end{equation}
where  $$D_{0} =\tau rI_{n_{2} }-\beta B^{T}B,\quad0.8\le\tau < 1\quad  and \quad r> \beta \left \| B^{T}B  \right \|.$$

He et al. proposed an optimal linearized ADMM in the literature \cite{22} with the iterative steps:
\begin{equation}
\left\{
\begin{aligned}
x^{k+1} &=arg\min_{x}\left \{ L_{\beta } \left ( x,y^{k},\lambda ^{k}   \right )|x\in X  \right \}
 \\y^{k+1} &=arg\min_{y}\left \{L_{\beta } \left ( x^{k+1} ,y,\lambda ^{k}   \right )+\frac{1}{2}\left \| y-y^{k}  \right \|_{D}^{2}   | y\in Y \right \}
 \\\lambda ^{k+1}& =\lambda ^{k} -\beta \left ( Ax^{k+1} +By^{k+1} -b \right ) 
\end{aligned}
\right.
\end{equation}
where$$D=\tau r I-\beta B^{T}B, \quad0.75< \tau < 1\quad and \quad\rho > \beta \left \| B^{T}B  \right \| $$
It is easy to see that the matrix $D$ is not necessarily semipositive definite. In the paper the authors prove that 0.75 is an optimal lower bound for $\tau$.

Another important issue for ADMM is to design an accelerated version of the original ADMM by slightly modifying it through a simple relaxation scheme. Note that ADMM is closely related to Proximity Point Algorithms\cite{23} (PPA). For the Proximity Point Algorithm, the convergence rate is improved if an additional over-relaxation step is added to the basic variables; see Tao\cite{24} for a relaxation variant of the PPA and proof of its linear convergence result.Gu and Yang\cite{25} further proved the optimal linear convergence rate of the relaxation PPA under regularity conditions. As Boyd et al. proved in\cite{3}, the execution of ADMM relies on the inputs of $\left ( y^{k},\lambda ^{k}   \right )$ and does not require $x^{k}$ at all. Thus, $x$ plays the role of an intermediate variable in \eqref{eqn:2}, while $\left ( y^{k},\lambda ^{k}   \right )$ is the basic variable. It is therefore natural to ask whether it is possible to include a relaxation step for the basic variables $\left ( y^{k},\lambda ^{k}   \right )$ in the ADMM scheme \eqref{eqn:2} to obtain a faster ADMM-type method. This idea leads to:
\begin{equation}
\left\{
\begin{aligned}
x^{k+1}&=arg\min_{x}\left \{ f\left ( x \right )-x^{T}A^{T}\lambda ^{k}+  \frac{\beta }{2}\left \| Ax+By^{k} -b \right \|^{2}        \right \}    
\\\hat{y}^{k}&=arg\min_{y}\left \{ g\left ( y \right )-y^{T}B^{T}\lambda ^{k}+  \frac{\beta }{2}\left \| Ax^{k+1}+By -b \right \|^{2}  \right \} 
 \\\hat{\lambda }^{k}&=\lambda ^{k}-\beta \left ( Ax^{k+1}+B\hat{y}^{k}-b  \right ) 
\end{aligned}
\right.
\end{equation}
\begin{equation}
\begin{pmatrix}y^{k+1} 
 \\\lambda ^{k+1} 
\end{pmatrix}=\begin{pmatrix}y^{k} 
 \\\lambda ^{k} 
\end{pmatrix}-\sigma \begin{pmatrix}y^{k}-\hat{y}^{k}  
 \\\lambda ^{k} -\hat{\lambda }^{k}  
\end{pmatrix}    
\end{equation}
where the relaxation factor $\sigma \in \left ( 0,2 \right ).$
Based on the above discussion,In this study, we introduce an adaptive linearized alternating direction multiplier method with a relaxation step that incorporates both a relaxation step and an adaptive technique. Our approach adopts a unified framework of variational inequalities and establishes the global convergence of this adaptive linearized alternating direction multiplier method with a relaxation step through theoretical analysis. Through the resolution of the Lasso problem, we demonstrate the algorithm's outstanding numerical performance. These findings will propel advancements in optimization algorithms and offer more efficient tools and methodologies for addressing practical challenges.

The rest of the paper is organized as follows. Section 2 gives some preliminaries; Section 3 gives iterative steps and proof of convergence of the algorithm; Section 4 gives the numerical experiments; and Section 5 draws the conclusions.

\section{Preliminaries}\label{sec2}
\noindent We present some preliminaries that will be used in convergence analysis.

\noindent In this article, the symbol $\left \| \cdot \right \|$ denotes the two-norm $\left \| \cdot \right \|_{2} .\left \| x \right \|_{D}^{2}: = x^{T}Dx$ is the matrix norm, where $D\in R^{n\times n} $ is a symmetric positive definite matrix and the vector $x\in R^{n}$. When $D $ is not a positive definite matrix, we still use the above notation.
\subsection{characterization of variational inequalities}\label{subsec2}

\noindent Define the Lagrangian function corresponding to problem \eqref{eqn:1} as:
\begin{equation}\label{eqn:14}
 L\left ( x ,y,\lambda   \right ) =f\left ( x \right )+g\left ( y \right ) -\lambda ^{T} \left ( Ax+By-b \right )
 \end{equation} 
In \eqref{eqn:14}, $(x,y)$ and $\lambda$ denote primitive and dual variables, respectively.

\noindent If
 $\left ( x^{\ast } ,y^{\ast },\lambda ^{\ast } \right )\in X\times Y\times R^{m} $ satisfies the following inequality:
$$ L\left ( x^{\ast } ,y^{\ast } ,\lambda  \right )\le L\left ( x^{\ast } ,y^{\ast } ,\lambda ^{\ast }  \right )\le L\left ( x ,y ,\lambda^{\ast }   \right ),\quad\forall \left ( x,y,\lambda  \right )\in X\times Y\times R^{m}=\Omega .$$
then $\left ( x^{\ast } ,y^{\ast } ,\lambda ^{\ast }  \right )$ is called a saddle point of  $L\left ( x ,y,\lambda   \right )$.\\
The above inequality is equivalent to the following variational inequality form:
\begin{equation}
\begin{cases}
x^{\ast } \in X,\quad f\left ( x \right )-  f\left ( x^{\ast }  \right )+\left ( x-x^{\ast }  \right )\left ( -A^{T} \lambda ^{\ast }  \right )\ge 0,\quad\forall x\in X . 
\\y^{\ast } \in Y,\quad g\left ( y \right )-  g\left ( y^{\ast }  \right )+\left ( y-y^{\ast }  \right )\left ( -A^{T} \lambda ^{\ast }  \right )\ge 0,\quad\forall y\in Y .
\\\lambda ^{\ast } \in R^{m} ,\qquad\left ( \lambda -\lambda ^{\ast }  \right )^{T}\left ( Ax^{\ast } +  By^{\ast } -b \right ) \ge 0,\quad\forall \lambda \in R^{m} .  
\end{cases}
\end{equation}
The above variational inequality can be written in the following form:
\begin{equation}\label{eqn:16}
VI\left ( \Omega ,F,\theta  \right ): \quad w^{\ast } \in \Omega ,\quad \theta \left ( u \right ) -\theta \left ( u^{\ast }  \right )+  \left ( w-w^{\ast }  \right ) ^{T}F\left ( w^{\ast }  \right ) \ge 0,\quad \forall w\in \Omega .   
\end{equation}
where\begin{equation}
\theta \left ( u \right ) =f\left ( x \right ) +  g\left ( y \right ) ,
u=\begin{pmatrix}x
 \\y
\end{pmatrix},v=\begin{pmatrix}y
 \\\lambda 

\end{pmatrix}, w=\begin{pmatrix}x
 \\y
 \\\lambda 
\end{pmatrix},F\left ( w \right ) =\begin{pmatrix}-A^{T}\lambda  
 \\-B^{T}\lambda  
 \\Ax+By-b
\end{pmatrix}    
\end{equation}
Owing to:$$\left ( w_{1}-w_{2}   \right )^{T}\left ( F\left (w_{1}  \right )-F\left ( w_{2}  \right )   \right )\ge 0,\quad\forall w_{1},w_{2}\in \Omega . $$
So, $F$ is a monotone operator.\\
The problem \eqref{eqn:1} is reformulated as the variational inequality \eqref{eqn:16} in this manner. We denote the set of solutions of the variational inequality $VI\left ( \Omega ,F,\theta  \right )$  as $\Omega ^{\ast } $ , where $\Omega ^{\ast }$ represents the set of non-empty solutions.

\subsection{some notation}\label{subsec2}
\noindent For the convenience of the proof, first define some auxiliary variables and matrices.
Let
\begin{equation}\label{eqn:18}
Q_{k+1}=\begin{pmatrix}
   \tau _{k}r\beta I_{n_{2} }  & O\\
  -B&\frac{1}{\beta }I_{m }  
\end{pmatrix}  \qquad and \qquad M=\begin{pmatrix}
  \sigma I_{n_{2} } &O \\
  -\sigma \beta B&\sigma I_{m }
\end{pmatrix}.
\end{equation}
\begin{equation}\label{eqn:19}
H_{k+1}=\frac{1}{\sigma } \begin{pmatrix}
  \tau_{k}r \beta I_{n_{2} }  & O\\
  O&\frac{1}{\beta } I_{m } 
\end{pmatrix} .
\end{equation}
\begin{equation}\label{eqn:120}
\tilde{w}^{k}=\begin{pmatrix}\tilde{x}^{k} 
 \\\tilde{y}^{k}
 \\\tilde{\lambda}^{k}
\end{pmatrix} =\begin{pmatrix}x^{k+1} 
 \\\hat{y}^{k} 
 \\\lambda  ^{k}  -\beta \left ( Ax^{k+1} +By ^{k} -b\right ) 
\end{pmatrix}.
\end{equation}
\begin{equation}
\tilde{v}^{k}=\begin{pmatrix}
 \tilde{y}^{k}
 \\\tilde{\lambda}^{k}
\end{pmatrix} =\begin{pmatrix} 
 \hat{y}^{k} 
 \\\lambda  ^{k}  -\beta \left ( Ax^{k+1} +By ^{k} -b\right ) 
\end{pmatrix}       
\end{equation}
\textbf{Lemma 2.1}. The $Q_{k+1},H_{k+1}$ and $M$ defined in \eqref{eqn:18} and \eqref{eqn:19} satisfy:
\begin{equation}\label{eqn:20}
 Q_{k+1}=H_{k+1}M\quad and\quad  H_{k+1}\succeq0,   
\end{equation}
Proof: \eqref{eqn:20} clearly hold. \\
\textbf{Lemma 2.2}. (Robbins-Siegmund Theorem\cite{26}) $a^{k}, b^{k}, c^{k} $ and $d^{k}$ are non-negative sequences and there are :
\begin{equation}
a^{k+1}\le \left ( 1+b^{k}  \right ) a^{k}+c^{k}-d^{k} ,\forall k=0,1,2\dots    
\end{equation}
if  $\sum_{k=0}^{+  \infty }b^{k} < +  \infty$ and $\sum_{k=0}^{+  \infty }c^{k} < +  \infty,$ so $\lim_{k \to \infty} a^{k}$ exists and is bounded, while there are $\sum_{k=0}^{+  \infty }d^{k} < +  \infty.$\\
\textbf{Definition 2.1}.If a function $f\left ( \cdot  \right )$ is a nonsmooth convex function on a convex set,if the following inequality holds:
$$f\left ( u \right ) \ge f\left ( v \right )+\phi ^{T} \left ( u-v \right ) ,\forall u\in \Omega .$$
then $\phi$  is said to be the subgradient of the function $f\left ( \cdot  \right )$ at $v\in \Omega $; The set consisting of all subgradients is called the subdifferential of the function $f\left ( \cdot  \right )$ at $v\in \Omega $, denoted $\partial{f\left ( v \right ) }$.

\subsection{Stopping criterion}\label{subsec3}
\noindent In the context of the two-block divisible convex optimization problem \eqref{eqn:1}, we investigate the optimality conditions for each subproblem during the iterative process of the Alternating Direction Method of Multipliers .\\
According to the optimality condition theorem, if $x^{\ast }$, $y^{\ast }$ are optimal solutions to the convex optimisation problem \eqref{eqn:1} and $\lambda ^{\ast }$ is the corresponding Lagrange multiplier, then the following conditions are satisfied.
\begin{equation}\label{eqn:24}
0\in \partial{f\left ( x^{\ast }  \right ) }-A^{T}\lambda ^{\ast} \end{equation}
\begin{equation}\label{eqn:25}
 0\in \partial{g\left ( y^{\ast }  \right ) }-B^{T}\lambda^{\ast} \end{equation}

\begin{equation}\label{eqn:26}
    Ax^{\ast } +By^{\ast } =b
\end{equation}
In this context, condition \eqref{eqn:26} is denoted as the original feasibility condition, while conditions \eqref{eqn:24} and \eqref{eqn:25} are labeled as the dual feasibility conditions.\\
According to the optimality conditions for the $y$-subproblem in \eqref{eqn:1}, we have:
\begin{equation}
\begin{aligned}
  0&\in \partial{g\left ( y^{k+1} \right ) }-B^{T}\lambda ^{k}+  \beta B^{T}\left ( Ax^{k+1}+By^{k+1}-b   \right )\\&= \partial{g\left ( y^{k+1} \right ) }-B^{T}\lambda ^{k+1}.       
\end{aligned}
\end{equation}
According to the optimality conditions for the $x$-subproblem in \eqref{eqn:1}, we have:
\begin{equation}
0\in \partial{f\left ( x^{k+1} \right ) }-A^{T}\lambda ^{k}+  \beta A^{T}\left ( Ax^{k+1}+By^{k}-b   \right ).    
\end{equation}
From the definition of $\lambda ^{k+1}$  in \eqref{eqn:1}, the above equation can be equated to:
\begin{equation}\label{eqn:29}
\begin{aligned}
 0&\in \partial{f\left ( x^{k+1} \right ) }-A^{T}\left ( \lambda ^{k}-\beta \left ( Ax^{k+1} +By^{k+1}-b \right )-\beta B\left ( y^{k} -y^{k+1}  \right )    \right )\\& = \partial{f\left ( x^{k+1} \right ) }-A^{T}\lambda ^{k+1}+\beta A^{T}B\left ( y^{k} -y^{k+1}  \right ) .   
\end{aligned}    
\end{equation}
\eqref{eqn:29} is equivalent to:
\begin{equation}\label{eqn:30}
  \beta A^{T}B\left ( y^{k+1}-y^{k}\right ) \in \partial{f\left ( x^{k+1} \right ) }-A^{T}\lambda ^{k+1}.  
\end{equation}
When comparing \eqref{eqn:30} with condition \eqref{eqn:24}, it is evident that the additional term is $\beta A^{T}B\left ( y^{k+1}-y^{k}\right ).$Thus, to verify dual feasibility, it is adequate to examine the residual $\beta A^{T}B\left ( y^{k+1}-y^{k}\right ).$

In summary, to ascertain the convergence of the alternating direction multiplier method, it is necessary to verify if the two residuals $p^{k+1}$  and $d^{k+1}$ are sufficiently small.\\
where$$p^{k+1} =\left \| Ax^{k+1}+By^{k+1}-b   \right \|.$$
$$d^{k} =\left \| \beta A^{T}B\left ( y^{k+1}-y^{k} \right )     \right \| .$$

 \section{Algorithm and convergence analysis}\label{sec3}
\subsection{New algorithms}\label{subsec2}
\noindent Suppose that $f:R^{n}\to R\cup \left \{ +\infty \right \}$ and $g:R^{n}\to R\cup \left \{ +\infty \right \}$ are proper, closed, convex functions. \\
The main iterative steps of the algorithm are:
\begin{spacing}{0.9}
\noindent\makebox[\textwidth]{\rule[0.1ex]{\textwidth}{0.01pt}}\\
Algorithm 1. An adaptive linearized alternating direction multiplier method with a relaxation step .\\[-1ex]
\makebox[\textwidth]{\rule[0.1ex]{\textwidth}{0.01pt}}
\end{spacing}
\noindent Set up: $\Omega  = X\times Y\times R^{m},\beta > 0,\tau _{k} > 0,p^{k}\ge 0 ,d^{k}\ge 0 ,r=\left \| B^{T}B  \right \|,\sigma \in \left ( 0,2 \right )  ,\varepsilon \in \left ( 0,2-\sigma  \right ) ,\Upsilon > 1,\rho > 1,\epsilon ^{pri} > 0,\epsilon ^{dual} > 0.$ choose  parameter sequence $\left \{ \eta _{k}  \right \}$ and $\left \{ s _{k}  \right \}$ , where $\sum_{k=0}^{\infty }\eta _{k}< +\infty$ and $\sum_{k=0}^{\infty }s _{k}< +\infty.$\\
\textbf{Step 0.} 
Input: $w^{0}=\left ( x^{0},y^{0},\lambda ^{0}    \right ) \in \Omega , \beta , \sigma, \tau _{0}, \tau _{min} , p^{0} , d^{0} , \Upsilon ,\varepsilon  , r, \rho , \epsilon ^{pri}  , \epsilon ^{dual} .$ set up $k=0. $\\
\textbf{Step 1.} Calculate: $w^{k+1}=\left ( x^{k+1},y^{k+1},\lambda ^{k+1}    \right ) \in X\times Y\times R^{m}.$  
\begin{equation}
\begin{cases}
x ^{k+1}=arg\min_{x}\left \{ f\left ( x\right ) -\left (\lambda ^{k}   \right )^{T}Ax  +  \frac{\beta }{2} \left \| Ax+By^{k} -b  \right \|^{2}  \right \}
 \\\hat{y}^{k}=arg\min_{y}\left \{ g\left ( y\right )-\left (\lambda ^{k}   \right )^{T}By+ \frac{\beta }{2}\left \| Ax^{k+1}+By-b  \right \|^{2}\right.\\\left.\qquad\qquad\qquad\qquad\qquad +\frac{1}{2}\left \| y-y^{k}  \right \|_{D_{k} }^{2}\right \} 
 \\\hat{ \lambda } ^{k}=\lambda ^{k}-\beta \left (  Ax ^{k+1}+B\hat{y}^{k}-b  \right ) 
\end{cases}
\end{equation}
where $D_{k}=\tau _{k}rI_{n_{2} }-\beta B^{T}B .$
\begin{equation}\label{eqn:32}
\begin{cases}
y^{k+1} =y^{k}-\sigma  \left (y^{k}-\hat{y}^{k}    \right ) \\
\lambda ^{k+1} =\lambda ^{k}-\sigma \left (\lambda ^{k}-\hat{\lambda }^{k}    \right ) 
\end{cases}
\end{equation}
where $\sigma  \in \left ( 0,2 \right )  .$\\
\textbf{Step 2.} If any one of the following conditions holds:
\begin{align}
&Condition 1.\Theta _{1}^{k} > \Theta _{2}^{k} . \\
&Condition 2.y^{k+1}=y^{k}. \nonumber
\end{align}
where $\Theta _{1}^{k} =\left ( 2-\sigma  \right ) \tau_{k}r\left \| y^{k} -y^{k+1}  \right \| ^{2}$, $  \Theta _{2}^{k}=\frac{1}{\varepsilon } \left \| B\left ( y^{k}-y^{k+1}  \right )  \right \|^{2}  .$\\
then go to step 3. otherwise $\tau _{k} =\gamma \ast\tau _{k}\left ( \gamma > 1 \right )  ,$ turn to step 1.\\
\textbf{Step 3.} If $\Theta _{1}^{k}- \Theta _{2}^{k}\ge \Upsilon \Theta _{2}^{k}.$ 
then set $ t_{k+1}=\max \left \{ \frac{r _{k} }{1+\eta _{k+1} } ,\tau _{min}  \right \}, $ Otherwise $t_{k+1}=\tau _{k} .$\\
\textbf{Step 4.} If any one of the following conditions holds: 
$$Condition 1.p^{k+1}> \left ( 1+s_{k}  \right )p^{k} .$$
$$Condition 2.d^{k+1}> \left ( 1+s_{k}  \right )d^{k} .$$
where $p^{k+1}=\left \| Ax^{k+1}+By^{k+1}-b   \right \|,d^{k+1}=\left \| \beta A^{T}B\left ( y^{k+1}-y^{k}   \right )   \right \|.$\\
then set $\tau _{k+1} =\rho \ast t _{k+1}\left ( \rho > 1 \right )  ,$ Otherwise $\tau_{k+1}=t _{k+1} .$\\ 
\textbf{Step 5.} If the stopping criterion is satisfied, return to step 6.\\
where, the stopping criterion is: 
\begin{equation}
\left \| Ax^{k+1}+By^{k+1}-b   \right \|\le \epsilon ^{pri}  \quad and \quad  \left \| \beta A^{T}B\left ( y^{k+1}-y^{k}   \right )   \right \| \le \epsilon ^{dual}     
\end{equation}
Otherwise make $k=k+1,$ and return to step 1.\\
\textbf{Step 6.} Output: $x^{k+1}  , y^{k+1},f\left ( x^{k+1}  \right ) +g\left ( y^{k+1}  \right )  .$

\begin{spacing}{0.9}
\noindent\makebox[\textwidth]{\rule[0.1ex]{\textwidth}{0.01pt}}\\
\end{spacing}
\subsection{Global convergence analysis}\label{subsec2}
\noindent In this section,we prove the global convergence for the proposed method.Before proceeding,we need the following lemma.\\
\textbf{Lemma 3.1}. Let the sequence $ \left \{ w^{k} \right \} $ be the iterative sequence generated by Algorithm 1, and $ \left \{ \tilde{w} ^{k} \right \} $ be defined in \eqref{eqn:120}, then, for any $ \tilde{w} ^{k} \in \Omega $, there are:
\begin{equation}\label{eqn:35}
\theta \left ( u \right ) -\theta \left ( \tilde{u} ^{k}  \right ) +\left ( w-\tilde{w}^{k}   \right )^{T} F\left ( \tilde{w}^{k} \right )  \ge \left ( v-\tilde{v}^{k}   \right )^{T}  Q_{k+1} \left (  v^{k}-\tilde{v}^{k}   \right ),\quad\forall w\in \Omega . 
\end{equation}
Proof: by the optimality condition for the $x$ subproblem in Algorithm 1: for any $ \tilde{x} ^{k}\in X $,there are:\\
$$f\left ( x \right )  -f\left ( \tilde{x} ^{k}  \right ) +\left ( x-\tilde{x} ^{k}  \right )^{T}\left \{ -A^{T}\lambda  ^{k}  +\beta A^{T}\left ( A\tilde{x}^{k} +By^{k} -b   \right ) \right \} \ge 0,\quad\forall x\in X.$$
Since $\tilde{\lambda }^{k}=\lambda ^{k}-\beta \left ( A\tilde{x}^{k} +By ^{k} -b\right ) $ in \eqref{eqn:120}, then: 
\begin{equation}\label{eqn:36}
 f\left ( x \right )  -f\left ( \tilde{x} ^{k}  \right ) +\left ( x-\tilde{x} ^{k}  \right )^{T}\left ( -A^{T}\tilde{\lambda } ^{k}   \right ) \ge 0.\quad \forall x\in X.    
\end{equation}
From the optimality conditions for the $y$-subproblem in Algorithm 1, it follows that for any $ \tilde{y} ^{k}\in Y $, there are:
$$g\left ( y \right ) -g\left (\tilde{y}^{k}  \right )+\left ( y-\tilde{y}^{k}    \right )^{T}\left \{-B^{T}\lambda ^{k}+  \beta B^{T}\left ( A\tilde{x}^{k} +  B\tilde{y}^{k}  -b   \right )\right.$$ $$\left.+  \left (   \tau_{k} rI_{n_{2} } -\beta B^{T}B\right )  \left ( \tilde{y}^{k}-y^{k}    \right )  \right \}\ge 0,\quad\forall y\in Y.$$
Also by the definition of $\tilde{\lambda }^{k}$, the above equation can be written as:
\begin{equation}\label{eqn:37}
 g\left ( y \right ) -g\left (\tilde{y}^{k}   \right )+\left ( y-\tilde{y}^{k}   \right )^{T}\left \{ -B^{T}\tilde{\lambda }^{k} +\tau_{k}r \beta    \left ( \tilde{y}^{k}-y ^{k}\right )\right \}\ge 0,\qquad\forall y\in Y.  \end{equation}
By the definition of $\tilde{\lambda }^{k}$ in \eqref{eqn:120}, we have:
$$\left ( A\tilde{x}+B\tilde{y}^{k}-b \right )-B\left ( \tilde{y}^{k}-y^{k}  \right )+  \frac{1}{\beta }\left ( \tilde{\lambda }^{k}-\lambda ^{k}    \right ) =0.$$
The above equation can be written as:
\begin{equation}\label{eqn：38}
\begin{split}
\tilde{\lambda } ^{k} \in R^{m},\qquad\left ( \lambda -\tilde{\lambda } ^{k} \right )^{T} \left \{ A\tilde{x}^{k}+B\tilde{y}^{k}-b-B\left ( \tilde{y}^{k}-y^{k}    \right )\right.\\\left.+  \frac{1}{\beta } \left ( \tilde{\lambda }^{k}-\lambda ^{k}    \right ) \right \} \ge 0,\quad\forall \lambda \in R^{m}. 
\end{split}
\end{equation}
The lemma can be proven by combining equations \eqref{eqn:36}, \eqref{eqn:37}, and \eqref{eqn：38} along with the notation $Q_{k+1}$.\\
\textbf{Lemma 3.2}. Let the sequence $\left \{ w^{k} \right \} $ be the iterative sequence generated by Algorithm 1 and $\left \{ \tilde{w} ^{k} \right \} $ be defined in \eqref{eqn:120}, then:
\begin{equation}\label{eqn:39}
    v^{k} -v^{k+  1} =M  \left ( v^{k}-\tilde{v}^{k}    \right )   
\end{equation}
where $M$ is defined in \eqref{eqn:18}. \\
Proof: This follows from  \eqref{eqn:32} and the definitions of $\tilde{\lambda }^{k}$ and $\hat{\lambda }^{k}$:\\
$\lambda ^{k+1} = \lambda ^{k}-\sigma \left ( \lambda ^{k}-\hat{\lambda}^{k}    \right ) $\\$ =\lambda ^{k} -\sigma  \beta \left ( A\tilde{x}^{k}+B\hat{y}^{k}-b   \right )$\\$= \lambda ^{k} -\sigma \left [\beta \left ( A\tilde{x}^{k}+By^{k}-b   \right ) -\beta B \left ( y^{k}-\hat{y}^{k}  \right )  \right ] $\\$=\lambda ^{k}-\sigma   \left ( \lambda ^{k}-\tilde{\lambda }^{k}    \right )+ \sigma   \beta B\left ( y^{k}-\tilde{y}^{k}    \right ) .$\\
This can be seen in combination with \eqref{eqn:32}:\\
$$\begin{pmatrix}y^{k+1} 
 \\\lambda ^{k+1} 

\end{pmatrix} = \begin{pmatrix}y^{k} 
 \\\lambda ^{k} 

\end{pmatrix}- \begin{pmatrix}
  \sigma  I_{n_{2} } &O \\
  -\sigma  \beta B&\sigma I_{m}
\end{pmatrix}  \begin{pmatrix}y^{k}-\tilde{y}^{k}   
 \\\lambda ^{k}-\tilde{\lambda }^{k}   
\end{pmatrix}.$$ 
Therefore, equation \eqref{eqn:39} holds.

It is clear from \eqref{eqn:39} and \eqref{eqn:20}:
\begin{equation*}
 \begin{aligned}
 \left ( v-\tilde{v}^{k}   \right )^{T} Q_{k+1} \left ( v^{k}- \tilde{v}^{k} \right )&= \left ( v-\tilde{v}^{k}   \right )^{T} H_{k+1} M\left ( v^{k}- \tilde{v}^{k} \right )\\&=\left ( v-\tilde{v}^{k}   \right )^{T} H_{k+1}\left ( v^{k}-v^{k+1}   \right ) .    
 \end{aligned}   
\end{equation*}
\textbf{Lemma 3.3}.Let the sequence $\left \{ w^{k} \right \} $ be the iterative sequence generated by Algorithm 1 and $\left \{ \tilde{w} ^{k} \right \} $ be defined in \eqref{eqn:120}, Then, we have:
\begin{equation}\label{eqn:40}
\begin{split}
   \left ( v-\tilde{v}^{k}   \right ) ^{T}H_{k+1}\left ( v^{k}-v^{k+1}   \right )&=\frac{1}{2} \left ( \left \| v-v^{k+1}  \right \|_{H_{k+1} }^{2}- \left \| v-v^{k}  \right \|_{H_{k+1} }^{2}  \right )\\&+  \frac{1}{2}\left \| v^{k}-\tilde{v}^{k}\right \|_{G_{k+1} }^{2},\qquad\forall v\in \Omega .
\end{split}
\end{equation}
where the matrix $G_{k+1}=Q_{k}^{T}+Q_{k}-M^{T}H_{k+1}M.$\\
Proof: Using the equation:
\begin{equation}\label{eqn:41}
\begin{split}
\left ( a-b \right )^{T}H_{k+1} \left (c-d\right ) &=\frac{1}{2}\left ( \left \| a-d \right \|_{H_{k+1}}^{2}-\left \| a-c \right \|_{H_{k+1}}^{2}   \right ) \\&+ \frac{1}{2}\left ( \left \| c-b \right \|_{H_{k+1}}^{2}-\left \| d-b \right \|_{H_{k+1}}^{2}   \right ).    
\end{split}
\end{equation}
In \eqref{eqn:41}, let $a=v$,$b=\tilde{v}^{k}$,$c=v^{k} $,$d=v^{k+1}$. We can get:
\begin{equation}\label{eqn:42}
\begin{split}
\left ( v-\tilde{v}^{k}\right )^{T}H_{k+1} \left (v^{k}-v^{k+1}\right )&=\frac{1}{2}\left ( \left \| v-v^{k+1} \right \|_{H_{k+1}}^{2}-\left \| v-v^{k} \right \|_{H_{k+1}}^{2}   \right ) \\&+ \frac{1}{2}\left ( \left \| v^{k}-\tilde{v}^{k} \right \|_{H_{k+1}}^{2}-\left \| v^{k+1}-\tilde{v}^{k} \right \|_{H_{k+1}}^{2}   \right ).  \end{split}
\end{equation}
For the last term in \eqref{eqn:42}, we have:
\begin{equation}\label{eqn:43}
\begin{split}
&\left \| v^{k}-\tilde{v}^{k} \right \|_{H_{k+1}}^{2}-\left \| v^{k+1}-\tilde{v}^{k} \right \|_{H_{k+1}}^{2}\\&=  \left \| v^{k}-\tilde{v}^{k} \right \|_{H_{k+1}}^{2}-\left \| \left ( v^{k} - \tilde{v}^{k} \right )-\left ( v^{k} -v^{k+1} \right )  \right \|_{H_{k+1}}^{2}\\& \stackrel{\left ( 44 \right ) }{=}\left \| v^{k}-\tilde{v}^{k} \right \|_{H_{k+1}}^{2}-\left \| \left ( v^{k} - \tilde{v}^{k} \right )-M\left ( v^{k} -\tilde{v}^{k}  \right )  \right \|_{H_{k+1}}^{2}\\&=2\left ( v^{k} -\tilde{v}^{k}   \right )^{T}H_{k+1}M\left ( v^{k} -\tilde{v}^{k}   \right ) - \left ( v^{k} -\tilde{v}^{k}   \right )^{T}M^{T} H_{k+1}M\left ( v^{k} -\tilde{v}^{k}   \right )\\&=\left ( v^{k} -\tilde{v}^{k}   \right )^{T}\left ( Q_{k}^{T}+Q_{k}-M^{T}H_{k+1}M \right ) \left ( v^{k} -\tilde{v}^{k}   \right ).   
\end{split}
\end{equation}
Combining \eqref{eqn:42} and \eqref{eqn:43}, \eqref{eqn:40} is proved. 

Now, let us examine the properties of the matrix $G_{k+1}$. Utilizing \eqref{eqn:20}, we can derive:
\begin{equation}
\begin{split}
 G_{k+1}&=Q_{k+1}^{T}+Q_{k+1}-M^{T}H_{k+1}M\\&=Q_{k+1}^{T}+Q_{k+1}-M^{T}Q_{k+1}\\&=\begin{pmatrix}
  2\tau_{k}r\beta I_{n_{2} }  & -B^{T} \\
  -B&\frac{2 }{\beta}I_{m}  
\end{pmatrix} -\begin{pmatrix}
  \sigma  I_{n_{2} }  & -\sigma \beta B^{T} \\
  O&\sigma  I_{m}  
\end{pmatrix}\begin{pmatrix}
 \tau_{k}r\beta I_{n_{2} }   & O \\
  -B&\frac{1}{\beta }  I_{m}  
\end{pmatrix}\\&=\begin{pmatrix}
  2\tau_{k}r\beta I_{n_{2} }  & -B^{T} \\
  -B&\frac{2}{\beta}I_{m}  
\end{pmatrix} - \begin{pmatrix}
  \sigma \tau_{k}r\beta I_{n_{2} }+\sigma \beta B^{T}B   & -\sigma B^{T} \\
 -\sigma B &\frac{\sigma }{\beta }I_{m}  
\end{pmatrix}\\&= \begin{pmatrix}
  \left ( 2-\sigma  \right )\tau_{k}r \beta I_{n_{2} }-\sigma \beta B^{T}B    & \left ( \sigma -1 \right )B^{T}  \\
  \left ( \sigma -1 \right )B&\frac{\left ( 2-\sigma  \right ) }{\beta }I_{m}  
\end{pmatrix}.
\end{split}    
\end{equation}

To simplify the proof, we give some properties of the parameters as follows:\\
From the iterative format of Algorithm 1:$\frac{1}{1+\eta _{k} }\tau _{k} \le \tau _{k+1} \le \left ( 1+\xi_{k} \right )\tau _{k} ,$That is, the sequence in Algorithm 1 satisfies:$\tau _{k}\subset \left [ \tau _{min},\tau _{max}   \right ].$ \\
Let the $(k+1)$th step ultimately satisfy the parameters of Step 4 as $\tau _{k+1} = \rho\gamma ^{m_{k}}\tau _{k}$, where $m_{k}$ is an integer. \\
Let $\xi _{k}: =\rho\gamma ^{m_{k} }-1,$ that is $1+\xi _{k}=\rho\gamma ^{m_{k} }.$ \\
then:$$\tau _{k+  1}\ge \frac{1+  \xi _{k} }{1+  \eta _{k} }\tau _{k}\ge \cdots\ge \frac{\prod_{i= 1}^{k} \left ( 1+  \xi _{i}  \right ) }{\prod_{i= 1}^{k} \left ( 1+  \eta  _{i}  \right )}\tau _{1} \ge \frac{\prod_{i= 1}^{k} \left ( 1+  \xi _{i}  \right ) }{\prod_{i= 1}^{\infty } \left ( 1+  \eta  _{i}  \right )}\tau _{0} .$$
By parameter setting $\sum_{k=1}^{\infty } \eta _{k} < +  \infty $  know $\prod_{i=1}^{\infty } \left ( 1+\eta _{i}  \right ) < +  \infty.$ \\
Making $k\to +  \infty$  gives: $\prod_{i=1}^{\infty } \left ( 1+  \xi _{i}  \right )< +  \infty$, that is: $ \sum_{k=1}^{\infty }    \xi _{i} < +  \infty.$\\
\textbf{Theorem 3.1}. Let the sequence $\left \{ w^{k} \right \} $ be the iterative sequence generated by Algorithm 1 and $\left \{ \tilde{w} ^{k} \right \} $ be defined in \eqref{eqn:120}. Then, we have:
\begin{equation}
\begin{split}
\left \| v^{k+1} -v^{\ast }  \right \|_{H_{k+1} }^{2} & \le \left ( 1+  \xi_{k}  \right )  \left \| v^{k} -v^{\ast }  \right \|_{H_{k} }^{2} -\frac{1}{\sigma ^{2} } \left \{ \beta \left \| y^{k}- y^{k+1} \right \|_{T_{k+1} }^{2}\right.\\&\left.+  \frac{\left ( 2-\sigma  \right )- \varepsilon }{\beta }\left \| \lambda ^{k}- \lambda ^{k+1} \right \|^{2}      \right \}   .  
\end{split}   
\end{equation}
where $\varepsilon \in \left ( 0,2-\sigma  \right )$ and $T_{k+1} \succ 0.$ \\ 
Proof: In the following, we will further investigate the $\left \| v^{k}-\tilde{v}^{k}\right \|_{G_{k+1} }^{2}$ term and show how it can be bounded.\\
Due to $G_{k+1}= \begin{pmatrix}
  \left ( 2-\sigma  \right )\tau_{k}r \beta I_{n_{2} }-\sigma \beta B^{T}B    & \left ( \sigma -1 \right )B^{T}  \\
  \left ( \sigma -1 \right )B&\frac{\left ( 2-\sigma  \right ) }{\beta }I_{m}  
\end{pmatrix}$ and $v=\begin{pmatrix}y
 \\\lambda 
\end{pmatrix},$ We have:
\begin{equation}
\begin{split}
\left \| v^{k}-\tilde{v}^{k}\right \|_{G_{k+1} }^{2} &= \left ( 2-\sigma  \right )\tau_{k}r\beta \left \| y^{k}-\tilde{y}^{k}\right \|^{2} -\sigma \beta \left \| B\left ( y^{k}-\tilde{y}^{k}\right )  \right \|^{2} \\
&+  \frac{2-\sigma }{\beta } \left \| \lambda ^{k}-\tilde{\lambda }^{k}    \right \|^{2} +  2\left ( \sigma -1 \right )\left ( \lambda ^{k}-\tilde{\lambda }^{k}    \right ) ^{T}B\left ( y^{k}-\tilde{y}^{k}\right )\\
&= \left ( 2-\sigma  \right )\tau_{k}r\beta \left \| y^{k}-\tilde{y}^{k}\right \|^{2}+  \left ( 2-\sigma  \right )\beta \left \| A\tilde{x}^{k}+B\tilde{y}^{k}-b\right \|^{2}\\&+2\beta \left ( A\tilde{x}^{k}+B\tilde{y}^{k}-b \right )^{T}B\left ( y^{k}-\tilde{y}^{k} \right )
\end{split}
\end{equation}
Also because $\hat{\lambda }^{k}=\lambda ^{k}-\beta \left ( A\tilde{x}^{k}+B\tilde{y}^{k}-b\right )$ and $y^{k}-\tilde{y}^{k} =\frac{1}{\sigma }\left ( y^{k}-y^{k+1}\right ) .$\\
 Therefore, we have:
\begin{equation}\label{eqn:47}
\begin{split}    
\left \| v^{k}-\tilde{v}^{k}\right \|_{G_{k+1} }^{2}&=\frac{\left ( 2-\sigma  \right )}{\sigma ^{2} } \tau_{k}r\beta \left \| y^{k}-y^{k+1}\right \|^{2}+  \frac{\left ( 2-\sigma  \right ) }{\beta }\left \| \lambda ^{k}-\hat{\lambda }^{k}\right \|^{2}\\&+ \frac{ 2}{\sigma } \left ( \lambda ^{k}-\hat{\lambda }^{k}\right )^{T}B\left ( y^{k}-y^{k+1}\right ).
\end{split}
\end{equation}
It is clear from \eqref{eqn:32}:
$$\lambda ^{k}-\hat{\lambda }^{k}= \frac{1}{\sigma } \left ( \lambda ^{k}-\lambda ^{k+1}\right )  .$$
Substituting the above equation into \eqref{eqn:47} gives:
\begin{equation}\label{eqn:48}
\begin{split} 
 \left \| v^{k}-\tilde{v}^{k}\right \|_{G_{k+1} }^{2}&=\frac{\left ( 2-\sigma  \right )}{\sigma ^{2} } \tau_{k}r\beta \left \| y^{k}-y^{k+1}\right \|^{2}+  \frac{\left ( 2-\sigma  \right ) }{\sigma ^{2} \beta }\left \| \lambda ^{k}-\lambda ^{k+1}\right \|^{2}\\&+ \frac{ 2}{\sigma^{2}  } \left ( \lambda ^{k}-\lambda ^{k+1}\right )^{T}B\left ( y^{k}-y^{k+1}\right ).
\end{split}
\end{equation}
In the following, we estimate $\left ( \lambda ^{k}-\lambda ^{k+1}\right )^{T}B\left ( y^{k}-y^{k+1}\right )$.\\
From the $Cauchy-Schwarz$ inequality: for $\forall \delta > 0$, there is :
\begin{equation}\label{eqn:49}
\left ( \lambda ^{k}-\lambda ^{k+1} \right ) ^{T} \left ( By^{k} - By^{k+1}\right )\ge -\frac{\delta }{\beta }\left \|  \lambda ^{k}-\lambda ^{k+1} \right \|^{2}-\frac{1}{4\delta }\beta \left \| B\left ( y^{k}-y^{k+1}   \right )  \right \|^{2} .
\end{equation}
Substituting \eqref{eqn:49} into \eqref{eqn:48} shows that:
\begin{equation}\label{eqn:50}
\begin{split}
\left \| v^{k}-\tilde{v}^{k}\right \|_{G_{k+1} }^{2}&\ge \frac{1}{\sigma ^{2} }\left [ \left ( 2-\sigma  \right )\tau_{k}r\beta \left \| y^{k}-y^{k+1}\right \|^{2}-\frac{1}{2\delta }\beta \left \| B\left ( y^{k}-y^{k+1}\right )  \right \|^{2}\right.\\&\left.+  \frac{\left ( 2-\sigma  \right )-2\delta  }{\beta }\left \| \lambda ^{k} - \lambda ^{k+  1} \right \|^{2} \right ]. 
\end{split}
\end{equation}
Owing to:
\begin{equation}\label{eqn:51}
 \left ( w_{1}-w_{2}   \right )^{T}\left ( F\left (w_{1}  \right )-F\left ( w_{2}  \right )   \right )\ge 0,\quad\forall w_{1},w_{2}\in \Omega .    
\end{equation}
It follows from \eqref{eqn:51} and \eqref{eqn:35}:
$$\theta \left ( u \right ) -\theta \left ( \tilde{u} ^{k}  \right ) +\left ( w-\tilde{w}^{k}   \right )^{T} F\left ( w \right )  \ge \left ( v-\tilde{v}^{k}   \right )^{T}  Q_{k+1} \left (  v^{k}-\tilde{v}^{k}   \right ),\quad\forall w\in \Omega . $$
Let the above equation $w = w^{\ast } $ and combine with \eqref{eqn:16} to show that:
\begin{equation}\label{eqn:52}
\left ( \tilde{v}^{k}-v^{\ast }  \right )^{T}  Q_{k+1} \left (  v^{k}-\tilde{v}^{k}   \right )\ge 0,\quad\forall w\in \Omega .  \end{equation}
 It is clear from \eqref{eqn:20} and \eqref{eqn:39}:
 \begin{equation}\label{eqn:53}
\left ( \tilde{v}^{k}-v^{\ast }  \right )^{T}  Q_{k+1} \left (  v^{k}-\tilde{v}^{k}   \right )=\left ( \tilde{v}^{k}-v^{\ast }  \right )^{T}  H_{k+1} \left (  v^{k}-v^{k+1}   \right ).
\end{equation}
 It is clear from \eqref{eqn:40}:
 \begin{equation}\label{eqn:54}
 \begin{split}
\left ( \tilde{v}^{k}-v^{\ast }  \right )^{T}  H_{k+1} \left (  v^{k}-v^{k+1}   \right )&=\frac{1}{2}\left ( \left \| v^{k}-v^{\ast }   \right \|_{H_{k+1} } ^{2}- \left \| v^{k+1}-v^{\ast }   \right \|_{H_{k+1} }^{2}  \right )\\& - \frac{1}{2  }\left \| v^{k}-\tilde{v}^{k}    \right \|_{G_{k+1} }^{2}.     
 \end{split}
\end{equation}
It follows from \eqref{eqn:52}, \eqref{eqn:53} and \eqref{eqn:54}:
\begin{equation}\label{eqn:55}
\left \| v^{k+1} -v^{\ast }  \right \|_{H_{k+1} }^{2} \le \left \| v^{k} -v^{\ast }  \right \|_{H_{k+1} }^{2}- \left \| v^{k} -\tilde{v}^{k}    \right \|_{G_{k+1} }^{2}.    
\end{equation}
Substituting \eqref{eqn:50} into \eqref{eqn:55} gives:
\begin{equation}\label{eqn:56}
\begin{split}
\left \| v^{k+1} -v^{\ast }  \right \|_{H_{k+1} }^{2} &\le \left \| v^{k} -v^{\ast }  \right \|_{H_{k+1} }^{2}-\frac{1}{\sigma ^{2} } \left \{ \beta \left \| y^{k}- y^{k+1} \right \|_{T_{k+1} }^{2}\right.\\&\left.+  \frac{\left ( 2-\sigma  \right )-2\delta  }{\beta }\left \| \lambda ^{k}- \lambda ^{k+1} \right \|^{2}      \right \} \\&\le \left \| v^{k} -v^{\ast }  \right \|_{H_{k} }^{2}+  \left ( \tau_{k}-\tau_{k-1}  \right ) \frac{r\beta }{\sigma }\left \| y^{k}-y^{\ast }   \right \|^{2}\\&-\frac{1}{\sigma ^{2} } \left \{ \beta \left \| y^{k}- y^{k+1} \right \|_{T_{k+1} }^{2}+  \frac{\left ( 2-\sigma  \right )-2\delta  }{\beta }\left \| \lambda ^{k}- \lambda ^{k+1} \right \|^{2}      \right \}\\& \le \left ( 1+  \xi_{k}  \right )  \left \| v^{k} -v^{\ast }  \right \|_{H_{k} }^{2} -\frac{1}{\sigma ^{2} } \left \{ \beta \left \| y^{k}- y^{k+1} \right \|_{T_{k+1} }^{2}\right.\\&\left.+  \frac{\left ( 2-\sigma  \right )-2\delta  }{\beta }\left \| \lambda ^{k}- \lambda ^{k+1} \right \|^{2}      \right \}  
\end{split} 
\end{equation}
where $T_{k+1} =\left ( 2-\sigma   \right ) \tau_{k}rI_{n_{2} }-\frac{1}{2\delta } B^{T} B  .$\\
For $\lim_{k \to \infty}\left ( v^{k}-v^{k+1}   \right ) =0 $ to hold, it is sufficient that:$\delta =\frac{\varepsilon }{2}$ and $T_{k+1}\succ 0.$\\
$T_{k+1}\succ 0$ is equivalent to $\left ( 2-\sigma   \right ) r_{k}\left \| y^{k}-y^{k+1}   \right \|^{2} > \frac{1}{\varepsilon  } \left \|  B \left ( y^{k}-y^{k+1} \right )  \right \| ^{2}.$\\
So \eqref{eqn:56} can be written as:
\begin{equation*}
\begin{split}
\left \| v^{k+1} -v^{\ast }  \right \|_{H_{k+1} }^{2} & \le \left ( 1+  \eta _{k}  \right )  \left \| v^{k} -v^{\ast }  \right \|_{H_{k} }^{2} -\frac{1}{\sigma ^{2} } \left \{ \beta \left \| y^{k}- y^{k+1} \right \|_{D_{k+1} }^{2}\right.\\&\left.+  \frac{\left ( 2-\sigma  \right )- \varepsilon }{\beta }\left \| \lambda ^{k}- \lambda ^{k+1} \right \|^{2}      \right \}   .  
\end{split}   
\end{equation*}
Therefore, Theorem 3.1 is proved.\\
\textbf{Theorem 3.2}. Let the sequence $\left \{ w^{k}  \right \} $ be the iterative sequence generated by Algorithm 1. Taking any point $w^{\ast } $ in $\Omega ^{\ast }$, we have: \\
(a).$ \lim_{ k\to \infty}\left \| v^{k+1} -v^{\ast }  \right \|_{H_{k+1} }^{2}$ exists and $\lim_{ k\to \infty}  \left \| v^{k+1} -v^{\ast }  \right \|_{H_{k+1} }^{2} < +  \infty.$ \\ 
(b).$\lim_{ k\to \infty}\left \| y^{k+1}-y^{k}   \right \|= 0$ and $\lim_{ k\to \infty}\left \| \lambda ^{k+1}-\lambda ^{k}   \right \|= 0 .$  \\
Proof: Let $$a^{k}= \left \| v^{k+1} -v^{\ast }  \right \|_{H_{k+1} }^{2},b^{k} = \xi _{k},c^{k} =0 ,$$$$d^{k}=\beta \left \| y^{k+1}-y^{k }    \right \| _{D_{k+1} }^{2} +  \frac{1-2\varepsilon }{\beta }\left \| \lambda ^{k+1}-\lambda ^{k }    \right \|^{2}.$$\\
By Lemma 2.2, (a) and (b) hold. \\
\textbf{Theorem 3.3}. Let the sequence $\left \{ w^{k}  \right \} $ be the iterative sequence generated by Algorithm 1, then $\left \{ w^{k}  \right \} $ converges to a point $w^{\infty } \in
\Omega ^{\ast }$.\\
Proof: It follows from Theorem 3.2: $\lim_{k \to \infty}\left \| v^{k+1}-v^{k}   \right \|=0.$\\ It follows from combining \eqref{eqn:39} and the non-singularity of $M$:$\lim_{k \to \infty} \left \| v^{k}-\tilde{v}^{k}    \right \|=0 .$\\
Since the sequence $\left \{ \left \| v^{k+1} -v^{\ast }  \right \|_{H_{k+1} }^{2} \right \}$ is a bounded sequence, \\ then for any fixed $v^{\ast }\in \Omega ^{\ast } ,$ there is a $\left \| v^{k+1}-v^{\ast }   \right \|$ bounded.\\
That is, the sequence $\left \{ v^{k} \right \}$ is bounded. \\
because of $$\left \| \tilde{v}^{k}-v^{\ast }    \right \|\le \left \| v^{k}-\tilde{v}^{k}    \right \|+\left \| v^{k}-v^{\ast }   \right \|.$$
It is known that $\left \| \tilde{v}^{k}-v^{\ast }    \right \|$ is bounded. \\
Clearly the sequence $\left \{ \tilde{v}^{k}   \right \} $ is also bounded and there must exist a convergence point $v^{\infty },\\$ such that a subsequence $\left \{ \tilde{v}^{k_{j} } \right \}$ of the existence sequence $\left \{ \tilde{v}^{k} \right \}$  converges to $v^{\infty } . $\\
From $\tilde{\lambda }^{k}=\lambda ^{k} -\beta \left ( A\tilde{x}^{k}+  By^{k}-b    \right ) $ in \eqref{eqn:120}, we have: \\$$A\tilde{x}^{k_{j} } =\frac{1}{\beta } \left ( \lambda ^{k_{j} } -\tilde{\lambda }^{k_{j} }   \right ) -\left ( By^{k_{j} }-b  \right ).$$\\
Since Matrix A is a column full rank matrix, it follows that the sequence $\left \{ \tilde{x}^{k_{j} }   \right \}$ converges. \\Set $ \lim_{k_{j}  \to \infty} \tilde{x}^{k_{j} } =x^{\infty } ,$ then there exists a subsequence $\left \{ \tilde{w}^{k_{j} }   \right \} $ converging to $w^{\infty }$.  \\  
Let $ k=k_{j} $ in \eqref{eqn:35}, we know that
$\qquad\tilde{w} ^{k_{j} }\in \Omega $,\\$$\theta \left ( u \right )-  \theta \left ( \tilde{u}^{k_{j} }  \right )+  \left ( w-\tilde{w}^{k_{j} }   \right )^{T}F\left ( \tilde{w}^{k_{j} }  \right )\ge \left ( w-\tilde{w}^{k_{j} }  \right ) ^{T}Q\left ( w^{k_{j} } -\tilde{w}^{k_{j} }  \right ) ,\quad\forall w\in \Omega . $$\\
Let the above equation $k \to \infty ,$ it is clear that:
\begin{equation}\label{eqn:57}
 w^{\infty } \in \Omega ,\theta \left ( u \right ) -\theta \left ( u^{\infty }  \right ) +  \left ( w-w^{\infty }  \right )^{T}F\left ( w^{\infty }  \right )\ge 0,\forall w\in \Omega .    
\end{equation}
\eqref{eqn:57} shows that $w^{\infty } $ in $\Omega ^{\ast }$ is a solution of the variational inequality $VI\left ( \Omega ,F,\theta \right )$, so the sequence $\left \{ w^{k} \right \} $ converges to $w^{\infty }. $

\section{Numerical experiments}\label{sec6}
\noindent The numerical performance of Algorithm 1 is presented in this paper through solving the Lasso problem and comparing Algorithm 1 with optimal linearized ADMM(OLADMM)\cite{22}. All simulation experiments were conducted on a laptop with 4GB of RAM memory using Matlab R2016a.

The LASSO problem\cite{12} in statistics is as follows:
\begin{equation}\label{eqn:58}
 \min_{x,y}\frac{1}{2}\left \| x-b \right \|_{2}^{2}   + \iota    \left \| y \right \|_{1} \qquad  s.t.\quad x=Ay   
\end{equation}
where $A\in R^{m\times n},b\in R^{m} ,x\in R^{n},y\in R^{n} .  $\\
Then applying Algorithm 1 to \eqref{eqn:58}, we obtain:

The $x$-subproblem is an unconstrained convex quadratic minimization problem, and its unique solution can be expressed as follows:
 $$x^{k+1} =\frac{1}{1+  \beta }\left ( b+  \lambda ^{k}+  \beta Ay^{k}   \right )$$

The $y$-subproblem is a $l_{1}+l_{2}$  minimization problem that can be solved using the soft-threshold operator:
\begin{equation*}
\begin{split}
  &\hat{y}^{k}   =shrink\left \{ y^{k}-\frac{1}{\tau _{k} \delta_{k}\beta   }A^{T} \left [ \lambda ^{k}-\beta \left ( x^{k+1}-Ay^{k}   \right )   \right ]  ,\frac{\iota }{\tau _{k} \delta_{k}\beta   }    \right \}.\\
 &y^{k+1} =y^{k}-\sigma \left ( y^{k}-\hat{y}^{k}    \right )  .    
\end{split}    
\end{equation*}

The solution to the $\lambda $-subproblem is:
\begin{equation*}
\begin{split}
&\hat{\lambda }^{k}=\lambda ^{k}-\beta \left ( x^{k+1} -A\hat{y}^{k}   \right )  .\\
&\lambda ^{k+1} =\lambda ^{k} -\sigma \left ( \lambda ^{k}-\hat{\lambda } ^{k}   \right ) .  
\end{split}
\end{equation*}

In this paper, numerical experiments are conducted to compare the performance of the  Algorithm 1 with OLADMM . \\\
Set the parameters of each algorithm as follows:\\\
Algorithm 1: $ r=\left \| A^{T}A \right \|$, $\iota=0.1\left \| A^{T}b  \right \|_{\infty }$,$\tau _{0} =0.75$,$\tau _{min} =0.01$,$\gamma  =1.2$  ,$\sigma =0.9$, $\beta =1$,$\rho = 3 $, $p^{0} =100$ , $d^{0} =100$ , $ {\varepsilon}' =\frac{1}{\varepsilon }=\frac{1}{2-\sigma }+  0.1  $, 
$\eta _{k}=0.25\min\left \{ 1,\frac{1}{\left ( \max\left \{ 1,k-l \right \}   \right )^{2}  }  \right \}  $ , $ s_{k}=2\min\left \{ 1,\frac{1}{\left ( \max\left \{ 1,k-l \right \}   \right )^{2}  }  \right \}.$\\
 where $ l $= the dimension of $\lambda $.\\
$\text{OLADMM: } r=\left \| A^{T}A  \right \|,\tau =0.75, \beta =1, \iota=0.1\left \| A^{T}b  \right \|_{\infty }$.\\
The stopping guidelines are: $$\left \| p^{k+1}  \right \| =\left \| x^{k+1} -Ay^{k+1}  \right \| < \varepsilon ^{pri} , \quad and \quad \left \| d^{k+1}  \right \|= \left \| \beta A \left ( y^{k+1}-y^{k}   \right )   \right \| < \varepsilon ^{dual} .$$
where $$\varepsilon^{pri}=\sqrt{n}\varepsilon ^{abs} +\varepsilon ^{rel} max\left \{ \left \| x^{k+1}  \right \| ,\left \| Ay^{k+1}  \right \|  \right \}$$ and $$\varepsilon ^{dual} =\sqrt{n}\varepsilon ^{abs} +\varepsilon ^{rel} \left \| y^{k+1}  \right \|.$$with $\varepsilon ^{abs}$ and $\varepsilon ^{rel}$ set to be $10^{-4} $ and$10^{-2} .$\\
For a matrix A of given dimension $m\times n$,we generate the data randomly as follow:\\$p=\frac{1}{n} ,$ $x^{0}=sprandn\left ( n,1,p \right ),$ $A=randn\left ( m,n \right )  , $ $b=A\ast x^{0}+sqrt\left ( 0.001 \right ) \ast  randn\left ( m,1 \right ) .$The initial
point is $\left ( y^{0},\lambda ^{0}   \right ) =\left ( 0,0 \right ) $. 

In order to observe the effect on the experiment caused by different values of $\sigma$ , the choice of the:$$\sigma =\left \{ 0.1,0.2,0.3,\cdots 1.7,1.8,1.9 \right \}$$ 
to carry out the experiment, the results of which show that when $\sigma \in \left [ 0.6,1.4 \right ]$, the results are satisfactory.In this paper, we take $\sigma =0.9.$

The experimental results of the two algorithms are shown below:\\
\text{Table 1.} Comparsion between Algorithm 1 and OLADMM for \eqref{eqn:58}. \\
\begin{tabular}{|c|c|c|c|c|c|c|c|c|c|}  
\hline \multicolumn{2}{|c|}{$n \times n$ matrix } & \multicolumn{4}{|c|}{ OLADMM } & \multicolumn{4}{|c|}{Algorithm 1} \\  
\hline$m$ & $n$ & Iter. & CPU(s) & $\left\|p^{k}\right\|$ & $\left\|q^{k}\right\|$ & Iter. & CPU(s) & $\left\|p^{k}\right\|$ & $\left\|q^{k}\right\|$ \\    
\hline 1000 & 1500 & 16 & 3.23 & $0.0700 $ & $0.0701$ & 11 & 2.38 & $0.0414$ & $0.0446$ \\  
\hline 1500 & 1500 & 13 & 2.81 & $0.0650$ & $0.0660$ & 10 & 2.53 & $0.0350$ & $0.0500$ \\  
\hline 1500 & 3000 & 17 & 7.03 & $0.0658$ & $0.0659$ & 10 & 6.18 & $0.0448$ & $0.0599$ \\  
\hline 2000 & 3000 & 14 & 7.16 & $0.0759$ & $0.0764$ & 10 & 6.57 & $0.0243$ & $0.0399$ \\  
\hline 3000 & 3000 & 13 & 7.80 & $0.0572$ & $0.0582$ & 9& 7.40 & $0.0377$ & $0.0695$ \\  
\hline 3000 & 5000 & 15 & 24.96 & $0.0685$ & $0.0688$ & 10& 23.44 & $0.0312$ & $0.0485$ \\
\hline 4000 & 5000 & 13 & 25.01 & $0.0664$ & $0.0674$ & 10& 24.69 & $0.0365$ & $0.0549$ \\
\hline 5000 & 5000 & 12 & 26.34 & $0.0614$ & $0.0636$ & 9& 25.73 & $0.0435$ & $0.0739$ \\
\hline
\end{tabular}

Table 1 presents the number of iterations and total computation time needed for OLADMM and Algorithm 1. It is evident that Algorithm 1 consistently outperforms OLADMM as it demands fewer iteration steps and less time to meet the termination condition. Consequently, Algorithm 1 demonstrates superior efficiency in contrast to OLADMM. Furthermore, Algorithm 1 significantly surpasses OLADMM, providing robust backing for our convergence analysis.

\section{Conclusion}\label{secA1}
In this study, we introduce an adaptive linearized alternating direction multiplier method with a relaxation step that incorporates both a relaxation step and an adaptive technique. Our approach adopts a unified framework of variational inequalities and establishes the global convergence of this adaptive linearized alternating direction multiplier method with a relaxation step through theoretical analysis. Through the resolution of the Lasso problem, we demonstrate the algorithm's outstanding numerical performance. These findings will propel advancements in optimization algorithms and offer more efficient tools and methodologies for addressing practical challenges.\\

\textbf{ Conflict of Interest}:The authors declare that they have no conflict of interest.






\bibliography{sn-bibliography}


\begin{thebibliography}{26}
\ifx \bisbn   \undefined \def \bisbn  #1{ISBN #1}\fi
\ifx \binits  \undefined \def \binits#1{#1}\fi
\ifx \bauthor  \undefined \def \bauthor#1{#1}\fi
\ifx \batitle  \undefined \def \batitle#1{#1}\fi
\ifx \bjtitle  \undefined \def \bjtitle#1{#1}\fi
\ifx \bvolume  \undefined \def \bvolume#1{\textbf{#1}}\fi
\ifx \byear  \undefined \def \byear#1{#1}\fi
\ifx \bissue  \undefined \def \bissue#1{#1}\fi
\ifx \bfpage  \undefined \def \bfpage#1{#1}\fi
\ifx \blpage  \undefined \def \blpage #1{#1}\fi
\ifx \burl  \undefined \def \burl#1{\textsf{#1}}\fi
\ifx \doiurl  \undefined \def \doiurl#1{\url{https://doi.org/#1}}\fi
\ifx \betal  \undefined \def \betal{\textit{et al.}}\fi
\ifx \binstitute  \undefined \def \binstitute#1{#1}\fi
\ifx \binstitutionaled  \undefined \def \binstitutionaled#1{#1}\fi
\ifx \bctitle  \undefined \def \bctitle#1{#1}\fi
\ifx \beditor  \undefined \def \beditor#1{#1}\fi
\ifx \bpublisher  \undefined \def \bpublisher#1{#1}\fi
\ifx \bbtitle  \undefined \def \bbtitle#1{#1}\fi
\ifx \bedition  \undefined \def \bedition#1{#1}\fi
\ifx \bseriesno  \undefined \def \bseriesno#1{#1}\fi
\ifx \blocation  \undefined \def \blocation#1{#1}\fi
\ifx \bsertitle  \undefined \def \bsertitle#1{#1}\fi
\ifx \bsnm \undefined \def \bsnm#1{#1}\fi
\ifx \bsuffix \undefined \def \bsuffix#1{#1}\fi
\ifx \bparticle \undefined \def \bparticle#1{#1}\fi
\ifx \barticle \undefined \def \barticle#1{#1}\fi
\bibcommenthead
\ifx \bconfdate \undefined \def \bconfdate #1{#1}\fi
\ifx \botherref \undefined \def \botherref #1{#1}\fi
\ifx \url \undefined \def \url#1{\textsf{#1}}\fi
\ifx \bchapter \undefined \def \bchapter#1{#1}\fi
\ifx \bbook \undefined \def \bbook#1{#1}\fi
\ifx \bcomment \undefined \def \bcomment#1{#1}\fi
\ifx \oauthor \undefined \def \oauthor#1{#1}\fi
\ifx \citeauthoryear \undefined \def \citeauthoryear#1{#1}\fi
\ifx \endbibitem  \undefined \def \endbibitem {}\fi
\ifx \bconflocation  \undefined \def \bconflocation#1{#1}\fi
\ifx \arxivurl  \undefined \def \arxivurl#1{\textsf{#1}}\fi
\csname PreBibitemsHook\endcsname

\bibitem[\protect\citeauthoryear{Tao et~al.}{2009}]{1}
\begin{botherref}
\oauthor{\bsnm{Tao}, \binits{M.}},
\oauthor{\bsnm{Yang}, \binits{J.}},
\oauthor{\bsnm{He}, \binits{B.}}:
Alternating direction algorithms for total variation deconvolution in image reconstruction.
TR0918, Department of Mathematics, Nanjing University
(2009)
\end{botherref}
\endbibitem

\bibitem[\protect\citeauthoryear{Ng et~al.}{2010}]{2}
\begin{barticle}
\bauthor{\bsnm{Ng}, \binits{M.K.}},
\bauthor{\bsnm{Weiss}, \binits{P.}},
\bauthor{\bsnm{Yuan}, \binits{X.}}:
\batitle{Solving constrained total-variation image restoration and reconstruction problems via alternating direction methods}.
\bjtitle{SIAM journal on Scientific Computing}
\bvolume{32}(\bissue{5}),
\bfpage{2710}--\blpage{2736}
(\byear{2010})
\end{barticle}
\endbibitem

\bibitem[\protect\citeauthoryear{Boyd et~al.}{2011}]{3}
\begin{barticle}
\bauthor{\bsnm{Boyd}, \binits{S.}},
\bauthor{\bsnm{Parikh}, \binits{N.}},
\bauthor{\bsnm{Chu}, \binits{E.}},
\bauthor{\bsnm{Peleato}, \binits{B.}},
\bauthor{\bsnm{Eckstein}, \binits{J.}}, \betal:
\batitle{Distributed optimization and statistical learning via the alternating direction method of multipliers}.
\bjtitle{Foundations and Trends{\textregistered} in Machine learning}
\bvolume{3}(\bissue{1}),
\bfpage{1}--\blpage{122}
(\byear{2011})
\end{barticle}
\endbibitem

\bibitem[\protect\citeauthoryear{Combettes and Pesquet}{2007}]{4}
\begin{barticle}
\bauthor{\bsnm{Combettes}, \binits{P.L.}},
\bauthor{\bsnm{Pesquet}, \binits{J.-C.}}:
\batitle{A douglas--rachford splitting approach to nonsmooth convex variational signal recovery}.
\bjtitle{IEEE Journal of Selected Topics in Signal Processing}
\bvolume{1}(\bissue{4}),
\bfpage{564}--\blpage{574}
(\byear{2007})
\end{barticle}
\endbibitem

\bibitem[\protect\citeauthoryear{Gabay and Mercier}{1976}]{5}
\begin{barticle}
\bauthor{\bsnm{Gabay}, \binits{D.}},
\bauthor{\bsnm{Mercier}, \binits{B.}}:
\batitle{A dual algorithm for the solution of nonlinear variational problems via finite element approximation}.
\bjtitle{Computers \& Mathematics with Applications}
\bvolume{2}(\bissue{1}),
\bfpage{17}--\blpage{40}
(\byear{1976})
\end{barticle}
\endbibitem

\bibitem[\protect\citeauthoryear{Glowinski and Marrocco}{1974}]{6}
\begin{barticle}
\bauthor{\bsnm{Glowinski}, \binits{R.}},
\bauthor{\bsnm{Marrocco}, \binits{A.}}:
\batitle{Analyse numerique du champ magnetique d'un alternateur par elements finis et sur-relaxation ponctuelle non lineaire}.
\bjtitle{Computer Methods in Applied Mechanics \& Engineering}
\bvolume{3}(\bissue{1}),
\bfpage{55}--\blpage{85}
(\byear{1974})
\end{barticle}
\endbibitem

\bibitem[\protect\citeauthoryear{He}{2018}]{7}
\begin{barticle}
\bauthor{\bsnm{He}, \binits{B.}}:
\batitle{My 20 years research on alternating directions method of multipliers}.
\bjtitle{Oper. Res. Trans}
\bvolume{22}(\bissue{1}),
\bfpage{1}--\blpage{31}
(\byear{2018})
\end{barticle}
\endbibitem

\bibitem[\protect\citeauthoryear{Eckstein and Yao}{2017}]{8}
\begin{barticle}
\bauthor{\bsnm{Eckstein}, \binits{J.}},
\bauthor{\bsnm{Yao}, \binits{W.}}:
\batitle{Approximate admm algorithms derived from lagrangian splitting}.
\bjtitle{Computational Optimization and Applications}
\bvolume{68},
\bfpage{363}--\blpage{405}
(\byear{2017})
\end{barticle}
\endbibitem

\bibitem[\protect\citeauthoryear{He et~al.}{2002}]{9}
\begin{barticle}
\bauthor{\bsnm{He}, \binits{B.}},
\bauthor{\bsnm{Liao}, \binits{L.-Z.}},
\bauthor{\bsnm{Han}, \binits{D.}},
\bauthor{\bsnm{Yang}, \binits{H.}}:
\batitle{A new inexact alternating directions method for monotone variational inequalities}.
\bjtitle{Mathematical Programming}
\bvolume{92},
\bfpage{103}--\blpage{118}
(\byear{2002})
\end{barticle}
\endbibitem

\bibitem[\protect\citeauthoryear{Ng et~al.}{2011}]{10}
\begin{barticle}
\bauthor{\bsnm{Ng}, \binits{M.K.}},
\bauthor{\bsnm{Wang}, \binits{F.}},
\bauthor{\bsnm{Yuan}, \binits{X.}}:
\batitle{Inexact alternating direction methods for image recovery}.
\bjtitle{SIAM Journal on Scientific Computing}
\bvolume{33}(\bissue{4}),
\bfpage{1643}--\blpage{1668}
(\byear{2011})
\end{barticle}
\endbibitem

\bibitem[\protect\citeauthoryear{Chan et~al.}{2012}]{11}
\begin{barticle}
\bauthor{\bsnm{Chan}, \binits{R.H.}},
\bauthor{\bsnm{Tao}, \binits{M.}},
\bauthor{\bsnm{Yuan}, \binits{X.}}:
\batitle{Linearized alternating direction method of multipliers for constrained linear least-squares problem}.
\bjtitle{East Asian Journal on Applied Mathematics}
\bvolume{2}(\bissue{4}),
\bfpage{326}--\blpage{341}
(\byear{2012})
\end{barticle}
\endbibitem

\bibitem[\protect\citeauthoryear{Tibshirani}{1996}]{12}
\begin{barticle}
\bauthor{\bsnm{Tibshirani}, \binits{R.}}:
\batitle{Regression shrinkage and selection via the lasso}.
\bjtitle{Journal of the Royal Statistical Society Series B: Statistical Methodology}
\bvolume{58}(\bissue{1}),
\bfpage{267}--\blpage{288}
(\byear{1996})
\end{barticle}
\endbibitem

\bibitem[\protect\citeauthoryear{Recht et~al.}{2010}]{13}
\begin{barticle}
\bauthor{\bsnm{Recht}, \binits{B.}},
\bauthor{\bsnm{Fazel}, \binits{M.}},
\bauthor{\bsnm{Parrilo}, \binits{P.A.}}:
\batitle{Guaranteed minimum-rank solutions of linear matrix equations via nuclear norm minimization}.
\bjtitle{SIAM review}
\bvolume{52}(\bissue{3}),
\bfpage{471}--\blpage{501}
(\byear{2010})
\end{barticle}
\endbibitem

\bibitem[\protect\citeauthoryear{Tao and Yuan}{2011}]{14}
\begin{barticle}
\bauthor{\bsnm{Tao}, \binits{M.}},
\bauthor{\bsnm{Yuan}, \binits{X.}}:
\batitle{Recovering low-rank and sparse components of matrices from incomplete and noisy observations}.
\bjtitle{SIAM Journal on Optimization}
\bvolume{21}(\bissue{1}),
\bfpage{57}--\blpage{81}
(\byear{2011})
\end{barticle}
\endbibitem

\bibitem[\protect\citeauthoryear{Fang et~al.}{2015}]{15}
\begin{barticle}
\bauthor{\bsnm{Fang}, \binits{E.X.}},
\bauthor{\bsnm{He}, \binits{B.}},
\bauthor{\bsnm{Liu}, \binits{H.}},
\bauthor{\bsnm{Yuan}, \binits{X.}}:
\batitle{Generalized alternating direction method of multipliers: new theoretical insights and applications}.
\bjtitle{Mathematical programming computation}
\bvolume{7}(\bissue{2}),
\bfpage{149}--\blpage{187}
(\byear{2015})
\end{barticle}
\endbibitem

\bibitem[\protect\citeauthoryear{Woo and Yun}{2013}]{16}
\begin{barticle}
\bauthor{\bsnm{Woo}, \binits{H.}},
\bauthor{\bsnm{Yun}, \binits{S.}}:
\batitle{Proximal linearized alternating direction method for multiplicative denoising}.
\bjtitle{SIAM Journal on Scientific Computing}
\bvolume{35}(\bissue{2}),
\bfpage{336}--\blpage{358}
(\byear{2013})
\end{barticle}
\endbibitem

\bibitem[\protect\citeauthoryear{He and Yuan}{2013}]{17}
\begin{barticle}
\bauthor{\bsnm{He}, \binits{B.}},
\bauthor{\bsnm{Yuan}, \binits{X.}}:
\batitle{Linearized alternating direction method of multipliers with gaussianback substitution for separable convex programming}.
\bjtitle{Numerical Algebra, Control and Optimization}
\bvolume{3}(\bissue{2}),
\bfpage{247}--\blpage{260}
(\byear{2013})
\end{barticle}
\endbibitem

\bibitem[\protect\citeauthoryear{Ouyang et~al.}{2015}]{18}
\begin{barticle}
\bauthor{\bsnm{Ouyang}, \binits{Y.}},
\bauthor{\bsnm{Chen}, \binits{Y.}},
\bauthor{\bsnm{Lan}, \binits{G.}},
\bauthor{\bsnm{Pasiliao~Jr}, \binits{E.}}:
\batitle{An accelerated linearized alternating direction method of multipliers}.
\bjtitle{SIAM Journal on Imaging Sciences}
\bvolume{8}(\bissue{1}),
\bfpage{644}--\blpage{681}
(\byear{2015})
\end{barticle}
\endbibitem

\bibitem[\protect\citeauthoryear{Chan et~al.}{2012}]{19}
\begin{barticle}
\bauthor{\bsnm{Chan}, \binits{R.H.}},
\bauthor{\bsnm{Tao}, \binits{M.}},
\bauthor{\bsnm{Yuan}, \binits{X.}}:
\batitle{Linearized alternating direction method of multipliers for constrained linear least-squares problem}.
\bjtitle{East Asian Journal on Applied Mathematics}
\bvolume{2}(\bissue{4}),
\bfpage{326}--\blpage{341}
(\byear{2012})
\end{barticle}
\endbibitem

\bibitem[\protect\citeauthoryear{Wang and Yuan}{2012}]{20}
\begin{barticle}
\bauthor{\bsnm{Wang}, \binits{X.}},
\bauthor{\bsnm{Yuan}, \binits{X.}}:
\batitle{The linearized alternating direction method of multipliers for dantzig selector}.
\bjtitle{SIAM Journal on Scientific Computing}
\bvolume{34}(\bissue{5}),
\bfpage{2792}--\blpage{2811}
(\byear{2012})
\end{barticle}
\endbibitem

\bibitem[\protect\citeauthoryear{He et~al.}{2016}]{21}
\begin{barticle}
\bauthor{\bsnm{He}, \binits{B.}},
\bauthor{\bsnm{Ma}, \binits{F.}},
\bauthor{\bsnm{Yuan}, \binits{X.}}:
\batitle{Linearized alternating direction method of multipliers via positive-indefinite proximal regularization for convex programming}.
\bjtitle{Avaliable on}
(\byear{2016})
\doiurl{https://optimization-online.org}
\end{barticle}
\endbibitem

\bibitem[\protect\citeauthoryear{He et~al.}{2020}]{22}
\begin{barticle}
\bauthor{\bsnm{He}, \binits{B.}},
\bauthor{\bsnm{Ma}, \binits{F.}},
\bauthor{\bsnm{Yuan}, \binits{X.}}:
\batitle{Optimally linearizing the alternating direction method of multipliers for convex programming}.
\bjtitle{Computational Optimization and Applications}
\bvolume{75}(\bissue{2}),
\bfpage{361}--\blpage{388}
(\byear{2020})
\end{barticle}
\endbibitem

\bibitem[\protect\citeauthoryear{He}{2015}]{23}
\begin{barticle}
\bauthor{\bsnm{He}, \binits{B.}}:
\batitle{Ppa-like contraction methods for convex optimization: a framework using variational inequality approach}.
\bjtitle{Journal of the Operations Research Society of China}
\bvolume{3},
\bfpage{391}--\blpage{420}
(\byear{2015})
\end{barticle}
\endbibitem

\bibitem[\protect\citeauthoryear{Tao and Yuan}{2018}]{24}
\begin{barticle}
\bauthor{\bsnm{Tao}, \binits{M.}},
\bauthor{\bsnm{Yuan}, \binits{X.}}:
\batitle{On the optimal linear convergence rate of a generalized proximal point algorithm}.
\bjtitle{Journal of Scientific Computing}
\bvolume{74},
\bfpage{826}--\blpage{850}
(\byear{2018})
\end{barticle}
\endbibitem

\bibitem[\protect\citeauthoryear{Gu and Yang}{2019}]{25}
\begin{botherref}
\oauthor{\bsnm{Gu}, \binits{G.}},
\oauthor{\bsnm{Yang}, \binits{J.}}:
On the optimal linear convergence factor of the relaxed proximal point algorithm for monotone inclusion problems.
arXiv preprint arXiv:1905.04537
(2019)
\end{botherref}
\endbibitem

\bibitem[\protect\citeauthoryear{Robbins~H}{1971}]{26}
\begin{botherref}
\oauthor{\bsnm{Robbins~H}, \binits{S.D.}}:
A convergence theorem for non negative almost supermartingales and some applications.
Optimizing Methods in Statistics,
233--257
(1971)
\end{botherref}
\endbibitem

\end{thebibliography}

\end{document}